\documentclass[12pt]{amsart}
\usepackage{amsfonts}

\theoremstyle{definition}

\theoremstyle{remark}

\begin{document}

\begin{center}

{\bf OSCILLATING OPERATORS IN BILATERAL GRAND LEBESGUE SPACES} \par

\vspace{3mm}
{\bf E. Ostrovsky}\\

e - mail: galo@list.ru \\

\vspace{3mm}

{\bf L. Sirota}\\

e - mail: sirota@zahav.net.il \\

\vspace{3mm}

 Abstract. \\
{\it In this paper we obtain the non - asymptotic estimations for oscillating
integral operators in the so - called
Bilateral Grand Lebesgue Spaces. We also give examples to show the
sharpness of these inequalities.} \\

\end{center}

\vspace{3mm}




\vspace{3mm}


\section{Introduction}

\vspace{3mm}

 The linear integral operator $ T_{\lambda} f(x), $ or, wore precisely, the
 {\it family of operators} of a view

 $$
 T_{\lambda}f(x) = \int_{R^d} \exp \left(i \lambda \Phi_1(x,y) \right) \
 \Phi_2(x,y) \ f(y) \ dy \eqno(0)
 $$
will be called {\it oscillating,} if $ \lambda $ is real "great" number: $ \lambda
> 2, \ \lambda \to \infty; \ \Phi_2(x,y) $ is a fixed non - zero smooth function:
$ \Phi_2(\cdot, \cdot) \in C_{\infty}^0  $ with finite support:

$$
\exists C = const \in (0,\infty), \ \Phi_2(x,y) = 0 \ \forall (x,y):
x^2 + y^2 \ge C, \eqno(1)
$$

 $ \Phi_1(x,y) $ is a fixed smooth function:
$ \Phi_1 \in C_{\infty}^0(R^d \times R^d) $ such that

$$
\det \left(  \frac{\partial^2 \Phi_1 }{ \partial x_i \ \partial y_j }  \right)
\ne 0 \eqno(2)
$$
on the support of the function $ \Phi_2(\cdot, \cdot). $ \par
 These operators are used in the theory of Fourier transform, theory of PDE,
 probability theory (study of characteristical functions and spectral densities) etc.\par
 In the physical applications the function $ \Phi_1 $ is called ordinary as {\it Phase
  function}, and the second function is named usually {\it Amplitude function.}\par
   The behavior of the function $ \lambda \to T_{\lambda} f $ as $ \lambda \to
   \infty  $ in the case when the function $ f(\cdot) $ is smooth, is described in
   the so - called stationary phase method (on the other words, saddle - point method). \par

 We denote as usually

 $$
 |f|_p = \left( \int_{R^d} |f(x)|^p \ dx  \right)^{1/p};  \  f \in L_p \
 \Leftrightarrow |f|_p < \infty.
 $$
We will consider further only the values $ p $ from the open interval
$ p \in (1,2) $ and denote $ q = q(p) = p/(p - 1); $ evidently,  $ q \in (2, \infty). $ \par

It is proved by E.M.Stein, see, e.g. in the book \cite{Stein1}, p. 307 - 355 that
the following estimation holds for the oscillating integral operator (0):

$$
|T_{\lambda} f|_{q} \le A(\Phi_1, \Phi_2)  \ \lambda^{-d/q} \ |f|_p, \eqno(3)
$$
$  A(\Phi_1, \Phi_2) \in (0, \infty). $  \par

 {\bf Our aim is a generalization of estimation (3) on the so - called Bilateral
 Grand Lebesgue Spaces $ BGL = BGL(\psi) = G(\psi), $ i.e. when } $ f(\cdot)
 \in G(\psi). \ $ \par
  We recall briefly the definition and needed properties of these spaces.
  More details see in the works \cite{Fiorenza1}, \cite{Fiorenza2}, \cite{Ivaniec1},  \cite{Ivaniec2}, \cite{Ostrovsky1}, \cite{Ostrovsky2}, \cite{Kozatchenko1},
  \cite{Jawerth1}, \cite{Karadzhov1} etc. More about rearrangement invariant spaces
  see in the monographs \cite{Bennet1}, \cite{Krein1}. \par

\vspace{3mm}

For $a$ and $b$ constants, $1 \le a < b \le \infty,$ let $\psi =
\psi(p),$ $p \in (a,b),$ be a continuous positive
function such that there exists a limits (finite or not)
$ \psi(a + 0)$ and $\psi(b-0),$  with conditions $ \inf_{p \in (a,b)} > 0 $ and
 $\min\{\psi(a+0), \psi(b-0)\}> 0.$  We will denote the set of all these functions
 as $ \Psi(a,b) = \Psi(a,b). $ \par

The Bilateral Grand Lebesgue Space (in notation BGLS) $  G(\psi; a,b) =
 G(\psi) $ is the space of all measurable
functions $ h: R^d \to{\mathbb R} $ endowed with the norm

$$
||h||G(\psi) \stackrel{def}{=}\sup_{p \in (a,b)}
\left[ \frac{ |h|_p}{\psi(p)} \right]. \eqno(4)
$$

The  $G(\psi)$ spaces with $ \mu(X) = 1$ appeared in \cite{Kozatchenko1}.
They are rearrangement invariant spaces and moreover interpolation spaces
between the spaces $ L_1(R^d) $ and $ L_{\infty}(R^d) $ under real interpolation
method \cite{Carro1}, \cite{Jawerth1}. \par
It was proved also that in this case each $ G(\psi) $ space coincides
with certain exponential Orlicz space, up to norm equivalence. In others
quoted publications were investigated, for instance,
 their associate spaces, fundamental functions
$\phi(G(\psi; a,b);\delta),$ Fourier and singular operators,
conditions for convergence and compactness, reflexivity and
separability, martingales in these spaces, etc.\par

{\bf Remark 1.} If we introduce the {\it discontinuous} function

$$
\psi_r(p) = 1, \ p = r; \psi_r(p) = \infty, \ p \ne r, \ p,r \in (a,b)
$$
and define formally  $ C/\infty = 0, \ C = const \in R^1, $ then  the norm
in the space $ G(\psi_r) $ coincides with the $ L_r $ norm:

$$
||f||G(\psi_r) = |f|_r.
$$

Thus, the Bilateral Grand Lebesgue spaces are direct generalization of the classical
Lebesgue spaces $ L_r. $ \par

 We recall the expression for the fundamental function for $ G(\psi) $ spaces.
Namely,

$$
\phi(G(\psi; a,b); \delta) = \sup_{p \in (a,b)}
\left[ \frac{\delta^{1/p}}{\psi(p)} \right]. \eqno(5)
$$
 More information about the fundamental function for $ G(\psi) $ spaces see in the
 article \cite{Ostrovsky2}; there was  considered, in particular,  many examples of
 $ G(\psi) $ spaces with exact calculation of their fundamental functions. \par

The BGLS norm estimates, in particular, Orlicz norm estimates for
measurable functions, e.g., for random variables are used in PDE
\cite{Fiorenza1}, \cite{Ivaniec1}, theory of probability in Banach spaces
\cite{Ledoux1}, \cite{Kozatchenko1},
\cite{Ostrovsky1}, in the modern non-parametrical statistics, for
example, in the so-called regression problem \cite{Ostrovsky1}.\par

The article is organized as follows. In the next section we obtain
the main result: upper bounds for oscillating operators in the Bilateral
Grand Lebesgue spaces. In the last section we  study the sharpness of the
obtained results by the building of the suitable examples.\par

 We use symbols $C(X,Y),$ $C(p,q;\psi),$ etc., to denote positive
constants along with parameters they depend on, or at least
dependence on which is essential in our study. To distinguish
between two different constants depending on the same parameters
we will additionally enumerate them, like $C_1(X,Y)$ and
$C_2(X,Y).$ The relation $ g(\cdot) \asymp h(\cdot), \ p \in (A,B), $
where $ g = g(p), \ h = h(p), \ g,h: (A,B) \to R_+, $
denotes as usually

$$
0< \inf_{p\in (A,B)} h(p)/g(p) \le\sup_{p \in(A,B)}h(p)/g(p)<\infty.
$$
The symbol $ \sim $ will denote usual equivalence in the limit
sense.

\bigskip

\section{Main result: upper estimations}

\vspace{3mm}

Let $ \psi(\cdot) \in \Psi(a,b), \ $ where $ 1 \le a < b \le 2. $ We define for the
values $ \lambda \ge 1 $  and the values $ q \in ( b/(b - 1), a/(a - 1) ), $ where
by definition at  $ a = 1 \ \Rightarrow a/(a - 1) = + \infty: $

$$
\psi^{(\lambda)}(q) = \lambda^{-d/q} \cdot \psi(q/(q-1)), \eqno(6)
$$
and define for the non - zero functions $ f $ belonging to the space $ G(\psi) $

$$
Z(\lambda, \psi, f) = \frac{||T_{\lambda} f ||G(\psi^{(\lambda)} ) }{ ||f||G(\psi) }. \eqno(7)
$$

{\bf Theorem 1.}

$$
\sup_{\lambda \ge 1} \sup_{\psi \in \Psi(a,b)}
 \sup_{f \in G(\psi), f \ne 0} Z(\lambda, \psi,f) \le A(\Phi_1, \Phi_2)
< \infty. \eqno(8).
$$

{\bf Proof.} Denote for the simplicity $ u = T_{\lambda}f; \ u: R^d \to R. $
 We can assume without loss of generality that $ ||f|| G(\psi) \le 1; $
this means that

$$
\forall p \in (a,b) \ \Rightarrow |f|_p \le \psi(p).
$$

Using the inequality (3), we obtain the estimation

$$
|u|_q \le A(\Phi_1, \Phi_2) \ \lambda^{- d/q } \ \psi(p) =
A(\Phi_1, \Phi_2) \ \lambda^{- d/q } \ \psi(q/(q - 1)) = A(\Phi_1, \Phi_2) \ \psi^{(\lambda)}(q). \eqno(9)
$$
The assertion of theorem 1 follows after dividing on the  $ \psi^{(\lambda)}(q), $
tacking the maximum on the $ q $ and on the basis of the definition of the
$ G(\psi) $ spaces.\hfill $\Box$

\vspace{3mm}

 Now we offer the another version of upper estimations for oscillating operator in the
Bilateral Grand Lebesgue spaces. Let $ \psi(\cdot), \nu(\cdot), \zeta(\cdot) $ be three
functions from the space $ \Psi(a,b), \ 1 \le a < b \le 2, $ such that

$$
\nu(p) = \psi(p) \cdot \zeta(p). \eqno(10)
$$
 Let us denote $ \nu^*(q) = \nu(q/(q - 1)), \ q \in (b/(b - 1), a/(a - 1)). $  \par

\newpage
 {\bf Theorem 2.}

 $$
 \lambda^d \ ||T_{\lambda} \ f||G(\nu^*) \le \phi(G(\zeta), \lambda^d) \cdot
 ||f||G(\psi). \eqno(11)
 $$
{\bf Proof.} We use again the Stein's estimation 3, which we rewrite as

$$
\lambda^d \cdot |u|_{q(p)} \le \lambda^{d/p} \cdot |f|_p \le \lambda^{d/p}
\cdot ||f||G(\psi) \cdot \psi(p),  \ p \in (a,b). \eqno(12)
$$
 We get after dividing both sides of inequality (12) on  the  $ \nu(p) $ and
 $ \lambda^{-d}: $

 $$
 \lambda^d \cdot \frac{|u|_q }{\nu(p)} \le ||f||G(\psi) \cdot
  \frac{\lambda^{d/p}}{ \zeta(p)}. \eqno(13)
 $$

 Tacking supremum of the bide sides of inequality (13) over the variable $ p; p \in
 (a,b), $ and tacking into account the definition of the fundamental function, we
 conclude:

 $$
 \lambda^d \cdot ||u||G(\nu(q/(q - 1) ) \le \phi \left(G(\zeta), \lambda^d \right) \
 \cdot ||f||G(\psi). \eqno(14)
 $$
 The last assertion (14) is equivalent to the proposition of theorem 2.
\hfill $\Box$

\bigskip

\section{Low bounds.}

\vspace{3mm}

In this section we built some examples in order to illustrate the exactness of upper
estimations.  It is sufficient to consider only the one - dimensional case: $ d = 1,$
i.e. $ x,y \in R^1. $ \par
We choose here the phase function $ \Phi_1 $ such that

$$
\Phi_1 = \Phi_1^{(0)} (x,y) = xy, \ (x,y) \in [-1,1]
$$
and

$$
\Phi_2 = \Phi_2^{(0)} (x,y) = 1, \ (x,y) \in [-1,1].
$$

Let us denote for the quoted values $ p, q(p), f \in L_p, f \ne 0 $
$$
W(\lambda, f,p) = \frac{|T_{\lambda} f |_q \ \lambda^{d/q}  }{ |f|_p }. \eqno(15)
$$
where as before $ p \in (1,2], \ q = q(p) = p/(p - 1) \in [2,\infty). $ \par
From the inequality of E.M.Stein (3) follows that

$$
\sup_{\lambda \ge 1} \sup_{p \in (1,2)} \sup_{f \in L_p, f \ne 0} W(\lambda, f,p)
\le A(\Phi_1^0, \Phi_2^0) < \infty. \eqno(16).
$$

We intend to prove an inverse inequality at the critical points $ \lambda \to
 \infty $ and $ p \to 2 - 0. $ \par

{\bf Theorem 3.}

$$
\underline{\lim}_{\lambda \to \infty} \underline{\lim}_{p \to 2 - 0}
 \sup_{f \in L_p, f \ne 0} W(\lambda, f,p)
\ge A_1(\Phi_1^0, \Phi_2^0) > 0. \eqno(17).
$$

{\bf Proof.}  Let us consider the function
$$
f(y) = f_0(y) = |y|^{-1/2}, \ |y| \in (0,1]
$$
and $ f_0(y) = 0 $ when $ y = 0 $ or $ |y| > 1. $ We have:

$$
|f_0|_p^p = 2 \int_0^1 y^{-p/2} \ dy = \frac{4}{2 - p}, p \in [1,2),
$$
or equally

$$
|f_0|_p =  [4/(2 - p)]^{1/p}. \eqno(18).
$$
 Further, let us investigate the function $ u = T_{\lambda} f. $
 Auxiliary denotation: $ \Lambda = \lambda x. $

$$
u = \int_{-1}^1 \exp(i \lambda x y) \ |y|^{-1/2} \ dy = 2 \int_0^1 \cos(\Lambda y) \
y^{-1/2} \ dy =
$$

$$
2 \Lambda^{-1/2} \ \int_0^{\Lambda} \cos z  /\sqrt{z} \ dz = 2 \Lambda^{-1/2} I(\Lambda),
$$
where
$$
I(\Lambda) = \int_0^{\Lambda} z^{-1/2} \ \cos z \ dz. \eqno(19)
$$

It is easy to calculate:

$$
I(\Lambda) \asymp \sqrt{\Lambda}, \ \Lambda \in (0,1); \
|I(\Lambda)| \asymp 1, \ \Lambda \in (0,1);
$$

therefore

$$
|u| \asymp 1, \ \Lambda \in (0,1); \ |u| \asymp \Lambda^{-1/2}, \ \Lambda \ge 1.
$$

Further,

$$
|u|_q^q \ge C^q \int_{1/\lambda}^{\infty} (\lambda x)^{-q/2} \ dx =
2 \ C^q \ \lambda^{-1} \ (q - 2)^{-1}.  \eqno(20)
$$

Substituting into the expression for the functional $ W, $ we get to the conclusion
of theorem 3 after simple computations. \par

 We can generalize the assertion of last assertion on the $ G(\psi) $ spaces
 as follows. Let us denote
$$
\psi_0(p) = [4/(2 - p)]^{1/p}, \ p \in (2, \infty), \eqno(21)
$$
then $ f_0(\cdot) \in G(\psi_0) $ and $ ||f_0||G(\psi_0) = 1. $ \par

{\bf Theorem 4.}

$$
 \underline{\lim}_{\lambda \to \infty}
 Z(\lambda, \psi_0, f_0) >  A_2(\Phi_1^0, \Phi_2^0) > 0.  \eqno(22).
$$

{\bf Proof.} From theorem 1 follows that $ u(\cdot) \in G(\psi_0^{(\lambda)} ) $
and $ ||u||G(\psi_0^{(\lambda)} ) = C_2 < \infty.  $ Let us now estimate the $ L_r $
norm of the function $ u $ from below. \par
 Since at $ x \ge 1/\lambda $

$$
u(x) \ge C (\lambda x)^{-1/2},
$$
we have for the values $ r > 2: $

$$
|u|_r^r \ge C \int_1^{\infty} (\lambda x)^{-r/2} \ dx = 2C \lambda^{-1} (r - 2)^{-1};
$$

$$
|u|_r \ge C \lambda^{-1/r} (r- 2)^{-1/r}. \eqno(23)
$$

Choosing the value $ r = q(p) = p/(p-1), $ we obtain on the basis of inequality (23):

$$
C_1^{-1} \lambda^{1/q} \ |u|_q /|f|_p \ge \frac{(2 - p)^{1/p} }{ (q - 2)^{1/q} } =
 (p-1)^{1/p -1},
$$
where $ C_1 $ does not depend on the $ p $ and $ \lambda. $ \par
 As long as

 $$
 \inf_{p \in (1,2) } (p - 1)^{1/p - 1}  = C_2 > 0,
 $$
we conclude

$$
\inf_{\lambda \ge 2} \inf_{p \in (1,2)} \lambda^{1/q} \ |u|_q /|f|_p \ge
C_2(\Phi_1^0, \Phi_2^0, \psi_0) > 0,
$$
or equally

$$
 \inf_{\lambda \ge 2} \frac{  ||T_{\lambda} f_0 ||G(\psi_0^{(\lambda)}) }{||f||G(\psi_0)} \ge C > 0,
$$
which is equivalent to the assertion of theorem 4.\par


\vspace{7mm}

\begin{thebibliography}{99}

\vspace{4mm}

\bibitem{Bennet1}
C. Bennet and R. Sharpley, {\it Interpolation of operators.}
Orlando, Academic Press Inc., 1988.
\bibitem{Carro1}
M. Carro and J. Martin, {\it Extrapolation theory for the real
interpolation method.} Collect. Math. {\bf 33}(2002), 163--186.
\bibitem{Fiorenza1}
A. Fiorenza. {\it Duality and reflexivity in grand Lebesgue
spaces.} Collect. Math. {\bf 51}(2000), 131--148.
\bibitem{Fiorenza2}
A. Fiorenza and G.E. Karadzhov, {\it Grand and small Lebesgue
spaces and their analogs.} Consiglio Nationale Delle Ricerche,
Instituto per le Applicazioni del Calcoto Mauro Picine", Sezione
di Napoli, Rapporto tecnico 272/03(2005).
\bibitem{Ivaniec1}
T. Iwaniec and C. Sbordone, {\it On the integrability of the
Jacobian under minimal hypotheses.} Arch. Rat.Mech. Anal., {\bf
119}(1992), 129--143.
\bibitem{Ivaniec2}
T. Iwaniec, P. Koskela and J. Onninen, {\it Mapping of Finite
Distortion: Monotonicity and Continuity.} Invent. Math. {\bf
144}(2001), 507--531.
\bibitem{Jawerth1}
B. Jawerth and M. Milman, {\it Extrapolation theory with
applications.} Mem. Amer. Math. Soc. {\bf 440}(1991).
\bibitem{Karadzov1}
G.E. Karadzhov and M. Milman, {\it Extrapolation theory: new
results and applications.} J. Approx. Theory, {\bf 113}(2005),
38--99.
\bibitem{Kozatchenko1}
Yu.V. Kozatchenko and E.I. Ostrovsky, {\it Banach spaces of random
variables of subgaussian type.} Theory Probab. Math. Stat., Kiev,
1985,  42--56 (Russian).
\bibitem{Krein1}
S.G. Krein, Yu. Petunin and E.M. Semenov, {\it Interpolation of
Linear operators.} New York, AMS, 1982.
\bibitem{Ledoux1}
M. Ledoux and M. Talagrand. {\it Probability in
 Banach Spaces.} Springer, Berlin, 1991.\\
\bibitem{Ostrovsky1}
E.I. Ostrovsky, {\it Exponential Estimations for Random Fields.}
Moscow - Obninsk, OINPE, 1999 (Russian).
\bibitem{Ostrovsky2}
E. Ostrovsky and L.Sirota, {\it Moment Banach spaces: theory and applications.}
HAIT Journal of Science and Engeneering, {\bf C}, Volume 4, Issues 1 - 2,
pp. 233 - 262, (2007).\\
\bibitem{Davis1}
  H.W. Davis, F.J.Murray, J.K.Weber. {\it Families of } $ L_p - $
{\it spaces with inductive and projective topologies.} Pacific J.Math.
v. 34, (1970), 619-638.
\bibitem{Jawerth1}
Jawerth B., Milman M.  {\it Extrapolation theory with applications.}
Mem. Amer. Math. Soc. 440 (1991).
\bibitem{Karadzhov1}
Karadzhov G.E., Milman M. {\it Extrapolation theory: new
results and applications. } J. Approx. Theory, 113 (2005), 38-99.
\bibitem{Stein1}
E.M.Stein. {\it Oscillating Integral in Fourier Analysis. } In: Beijing Lectures in Harmonic Analysis, Princeton University  Press, (1986), p. 307 - 355.\\


\end{thebibliography}
\end{document}